\documentstyle[12pt]{article}
\begin{document}
\title{\Large \bf The Center for  the Elliptic Quantum Group $E_{\tau,\eta}(sl_n)$}
\author{{Shaoyou Zhao\footnote{E-mail: zsy@biophysics.nju.edu.cn}}\\
{\small  P.O. Box 2324, Department of Physics, Nanjing University}\\
{\small Nanjing, 210093, P.R. China}\\[3mm]
{ Kangjie Shi, Ruihong Yue}\\
{\small P.O.Box 105, Institute of modern physics, Northwest university }\\
{\small Xi'an, 710069, P.R. China}         
}
\maketitle
\begin{abstract}
  We give the center of the elliptic quantum group in general case. Based on the Dynamic Yang-Baxter Relation and the fusion method, we prove that the center commute with all generators of the elliptic quantum group. Then for a kind of assumed form of these generators,  we find that the coefficients of these generators form a new type closed algebra.  We also give the center for the algebra.
\end{abstract}

\section {Introduction}
~~~~
 Recently, many papers  have  focused on the many-body long-distance
integrable dynamical system, such as the Ruijsenaar Schneider model and
the Calogero Moser (CM) model[1-3]. They are
closely connected with the quantum Hall effect in the condense matter
physics and the Seiberg Witten (SW) theory in the field theory, especially
for the equations of the spectral curve in the SW theory, namely, the
modified eigenvalue equations of the Lax matrices in the above integrable
models[4-6]. These Lax
matrices are the classical limit of the $L$-matrices which is associated
with the interaction-round-a-face (IRF) model of Lie
group[7-9] and the
modified Yang-Baxter relation (NSF equation)[10-13].  
All these $L$-matrices are corresponding to the 
representation of the elliptic quantum group which was proposed by Felder and Varchenko\cite{fe1,fe2}.  The elliptic quantum group is an algebraic structure underlying the elliptic solution of the Yang-Baxter relation in the statistical mechanics and connected withe the Knizhnik-Zamolodchikov-Bernard equation on torus.
 In the theory of the group and  algebra, The center of them play an important role.
The center is generally defined by the elements which can commute with all generators of the group and algebra.  In Ref. \cite{fe2}, Felder and Varchenko obtained the center (also called the quantum determinant) for $n=2$ case. In first part of this paper, we will discuss the center of the elliptic quantum group for the general case. And in the second part of this paper,  assuming that the generators of the elliptic quantum group( namely, the elements of the $L$-matrices) take a certain kind of form, we find that their coefficients can construct a closed algebra. To this algebra, we also obtain its center.

\section{The DYBR and the   elliptic quantum group $E_{\tau,\eta}(sl_n)$}   
 
~~
It is well known that the Boltzmann weight of the  $A^{(1)}_{n-1}$
IRF\cite{jim1,jim2,jim3} model can be written as 
\begin{equation}\begin{array}{ll}
 R(a|z)_{ii}^{ii}=\displaystyle\frac{\sigma(z+w)}{\sigma(w)},&
 R(a|z)_{ij}^{ij}=\displaystyle\frac{\sigma(z)\sigma(
                 a_{ij}-w)}{\sigma(w)\sigma(a_{ij})}\quad
            \mbox{for}\;\;i\ne j,\\
R(a|z)_{ij}^{ji}=\displaystyle\frac{\sigma(z+a_{ij})}{\sigma(a_{ij})}
           \quad  \mbox{for}\;\;i\ne j,&
R(a|z)_{i'j'}^{i\ j}=0 \quad\mbox{for other cases},
\end{array}\end{equation}
where $a\equiv(m_0,m_1,\cdots,m_{n-1})$ is an $n$-vector, and $
a_{ij}= a_i-a_j$, $a_i=w(m_i-\frac{1}{n}\sum_l m_l+w_i)$,  $m_i$
($i=0,1,\cdots,n-1$) are integers which describe the state of model, while
 $\{w,w_i\}$ are generic c-numbers which are the parameters of the model,
and $\sigma(z)\equiv\theta\left[^{\frac{1}{2}}_{\frac{1}{2}}\right](z,\tau),$
with
$$\theta\left[^a_b\right](z,\tau)\equiv
             	\sum_{m\in Z}e^{i\pi(m+a)^2\tau+2i\pi(m+a)(z+b)}.$$
We define  an n-dimension vector $\hat j=(0,0,\cdots,0,1,0,\cdots)$, in
which $j$th component is 1.

We consider a matrix whose elements are linear operators acting on the quantum space $V_0$. The elements of the matrix are denoted by $ L(^a_b|z)^j_i$. The $R$-matrix and the $L$-matrix can also be depicted by the figure 1 and figure 2, respectively, and they satisfy the Well-known Dynamic Yang-Baxter Relation (DYBR). The DYBR is written as (figure 3)
\begin{equation}
\sum_{i',j'}R(b|z_1-z_2)_{ij}^{i'j'}L(_b^a|z_1)_{i'}^{i''}L(_{b+\hat
               i'}^{a+\hat i''}|z_2)_{j'}^{j''}
 =\sum_{i',j'}L(_b^a|z_2)_{j}^{j'}L(_{b+\hat j}^{a+\hat j'}|z_1)_{i}^{i'}
                      R(a|z_1-z_2)_{i'j'}^{i''j''},
\end{equation}
where $b\equiv(m_0^b,m_1^b,\cdots,m_{n-1}^b),$
$a\equiv(m_0^a,m_1^a,\cdots,m_{n-1}^a),$ We note that Eq.(2) gives the
quadratic relation of the elements of $L$. If we let
$b=a+h$, the form of the equation will be the same as that given
in the Ref.\cite{fe1,fe2}.\\
 
\begin{picture}(50,50)(0,0)
\put(150,-20){\vector(-2,1){80}}
\put(150,20){\vector(-2,-1){80}}
\put(107,-10){$a$}\put(127,-20){$i'$}\put(87,16){$i$}
\put(127,16){$j'$}\put(87,-20){$j$}\put(65,25){$z_1$}\put(65,-26.5){$z_2$}
\end{picture}
\begin{picture}(50,50)(0,0)
\put(300,5){\vector(-1,0){80}}
\put(260,35){\vector(0,-1){60}}\put(260.5,35){\vector(0,-1){60}}
\put(261,35){\vector(0,-1){60}}
\put(230,20){$b+{\hat i}$}
\put(280,-10){$a$}\put(210,-15){$b\equiv a+{\hat h}$}\put(240,8){$i$}
\put(275,8){$j$}\put(265,-20){$h$}\put(213,2.5){$z$}
\put(270,25){$a+{\hat j}$}
\end{picture}
\\ 
\bigskip
\\
$$\begin{array}{cc}
Figure\; 1: The\;elements \;of\;R\!-\!matrix\;\;\;\;&\;\;\;\;
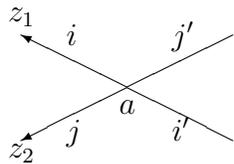
Figure\; 2: The\; element \;of\;L\!-\!matrix,\\
R(a|z_1-z_2)_{ij}^{i'j'}.&L(a,h|z)_i^j\equiv L(_b^a|z)_i^j.
\end{array}$$

$$
\begin{picture}(50,50)(0,0)
\put(-10,-32.5){\vector(-2,1){100}}
\put(-10,32.5){\vector(-2,-1){100}}
\put(-30,50){\vector(0,-1){100}}\put(-30.5,50){\vector(0,-1){100}}
\put(-31,50){\vector(0,-1){100}}
\put(-27,-35){$a$}\put(-95,-35){$b\equiv a+{\hat h}$}\put(-100,-21){$j$}
\put(-100,17){$i$}\put(-40,-40){$h$}\put(-120,19){$z_1$}
\put(-60,-5){$b+\hat i'$}
\put(-120,-19){$z_2$}\put(-20,-25){$i''$}\put(-25,30){$j''$}
\put(-23,0){$a+{\hat i''}$}\put(-65,10){$j'$}\put(-65,-18){$i'$}
\end{picture}
=
\begin{picture}(50,50)(0,0)
\put(150,-15.5){\vector(-2,1){100}}   
\put(150,15.5){\vector(-2,-1){100}}
\put(80,50){\vector(0,-1){100}}\put(80.5,50){\vector(0,-1){100}}
\put(81,50){\vector(0,-1){100}}
\put(116,-10){$a$}\put(30,-50){$b\equiv a+{\hat h}$}\put(60,18){$i$}
\put(38,0){$b+\hat j$}\put(82,-5){$a+\hat j'$}
\put(60,-24){$j$}\put(84,-40){$h$}\put(38,33){$z_1$}
\put(38,-37){$z_2$}\put(140,-25){$i''$}\put(140,17){$j''$}
\put(95,15){$i'$}\put(95,-25){$j'$}
\end{picture}$$
\\ \\
\\$$ Figure\; 3: The\;\; \;dynamical \;\;Yang-Baxter \;\;relation.$$

The elliptic quantum group $E_{\tau,\eta}(sl_n)$ is an algebra generated by the matrix elements of the $L$-matrix. For a given linear space $V$,  the generators of the algebra satisfy the following relations
\begin{eqnarray*}
& &R(b|z_1-z_2)^{ii}_{ii}
  L(_b^a|z_1)_{i}^{i''}L(_{b+\hat i}^{a+\hat i''}|z_2)_{i }^{i''}
=R(a|z_1-z_2)^{i''i''}_{i''i''}
  L(_b^a|z_2)_{i}^{i''}L(_{b+\hat i}^{a+\hat i''}|z_1)_{i}^{i''}\\
\nonumber\\
& &R(b|z_1-z_2)_{ii}^{ii}
  L(_b^a|z_1)_i^{i''}L(_{b+\hat i}^{a+\hat i''}|z_2)_i^{j''}\nonumber\\
&=&R(a|z_1-z_2)_{i''j''}^{i''j''}
  L(_b^a|z_2)_i^{j''}L(_{b+\hat i}^{a+\hat j''}|z_1)_i^{i''}
+R(a|z_1-z_2)_{j''i''}^{i''j''}
  L(_b^a|z_2)_i^{i''}L(_{b+\hat i}^{a+\hat i''}|z_1)_i^{j''}\nonumber\\ 
   & & \quad\quad\quad\quad\quad    (i''\ne j'') \\
& &R(a|z_1-z_2)_{i''i''}^{i''i''}
  L(_b^a|z_2)_j^{i''}L(_{b+\hat j}^{a+\hat i''}|z_1)_i^{i''}\nonumber\\
&=&R(b|z_1-z_2)_{ij}^{ij}
  L(_b^a|z_1)_i^{i''}L(_{b+\hat i}^{a+\hat i''}|z_2)_j^{i''}
+R(b|z_1-z_2)_{ij}^{ji}
  L(_b^a|z_1)_j^{i''}L(_{b+\hat j}^{a+\hat i''}|z_2)_i^{i''}\nonumber\\ 
&&\quad\quad\quad\quad\quad\quad     (i\ne j)  \\
& &R(b|z_1-z_2)_{ij}^{ij}
  L(_b^a|z_1)_i^{i''}L(_{b+\hat i}^{a+\hat i''}|z_2)_j^{j''}
+R(b|z_1-z_2)_{ij}^{ji}
  L(_b^a|z_1)_j^{i''}L(_{b+\hat j}^{a+\hat i''}|z_2)_i^{j''}\nonumber\\
&=&R(a|z_1-z_2)_{i''j''}^{i''j''}
  L(_b^a|z_2)_j^{j''}L(_{b+\hat j}^{a+\hat j''}|z_1)_i^{i''}
+R(a|z_1-z_2)_{j''i''}^{i''j''}
  L(_b^a|z_2)_j^{i''}L(_{b+\hat j}^{a+\hat i''}|z_1)_i^{j''}\nonumber\\ 
& & \;\;\;\;\;\;\;\;\;\;(i\ne j,\;i''\ne j'').
\end{eqnarray*}

\section{The center for the   elliptic quantum group}
~~
In this section, we will discuss the center of the  elliptic quantum group. For a given space $V\otimes V\otimes\cdots\otimes V\equiv V^{\otimes n}$ with a base $e_{i_1}\otimes e_{i_2}\otimes\cdots\otimes e_{i_n}$, we define the permutation operator $P$ as 
\begin{equation}
 P(^{1\ \ 2\ \cdots\ n}_{j_1\ j_2\ \cdots\ j_n})
   e_{i_1}\otimes e_{i_2}\otimes\cdots\otimes e_{i_n}
=e_{i_{j_1}}\otimes e_{i_{j_2}}\otimes\cdots\otimes e_{i_{j_n}}.
\end{equation}
Define 
\begin{equation}
P_-=\frac{1}{n!}\sum_P(-1)^{\delta_P}P,
\end{equation}
where $\delta_P=0$ if $P$ is even permutation, otherwise $\delta_P=1$.
Then, one can easily prove the $P_-$ satisfies the property of projection operator  
\begin{equation}
P_-P_-=P_-,
\end{equation}

Let $z=-w$ in the Eq.(1), we have
\begin{equation}\begin{array}{lll}
R(a|-w)^{ii}_{ii}=0, \quad &
\displaystyle R(a|-w)^{ij}_{ij}=-\frac{\sigma(a_{ij}-w)}{\sigma(a_{ij})},  \quad&
\displaystyle R(a|-w)^{ji}_{ij}=\frac{\sigma(a_{ij}-w)}{\sigma(a_{ij})},
\end{array}\end{equation}
where $i\ne j$,  the other elements of $R-$matrix are all zeros. So we can easily obtain 
\begin{equation}
R(a|-w)P=-R(a|-w).
\end{equation}
Then for the above space $V^{\otimes n}$, we define a Cherednik operator $A$ as
\begin{eqnarray}
A&=&R_{12}(a|-w)R_{13}(a|-2w)\cdots R_{1N}(a|-(n-1)w)\nonumber\\
 & &\times R_{23}(a|-w)R_{24}(a|-2w)\cdots R_{2N}(a|-(n-2)w)\nonumber\\
 & &\times \cdots R_{n-1,n}(a|-w).
\end{eqnarray}
Similarly with the Eq.(7), one can check that
\begin{equation}
AP_-=A.
\end{equation}
On the other hand, we can always find an invertible matrix $B$ satisfying $A=BP_-$.  Define ${\bf L}$ in the space $V_0\otimes V^{\otimes n}$ as in the figugre 4. Taking use of the DYBR repeatedly, we can find a relation between $A$ and ${\bf L}$, i.e., $A{\bf L}={\bf L}'A$, where $L'$ is defined as figure 5. Therefore, applying the projection operator $P_-$ on ${\bf L}$, we have
\begin{equation}
P_-{\bf L}=B^{-1}A{\bf L}=B^{-1}{\bf L}'A=B^{-1}{\bf L}'AP_-=P_-{\bf L}P_-.
\end{equation}
\begin{picture}(50,50)(0,0)
\put(40,5){\vector(1,0){120}}\put(40,5.3){\vector(1,0){120}}
\put(40,5.6){\vector(1,0){120}}\put(40,5.9){\vector(1,0){120}}
\put(125,35){\vector(0,-1){60}}
\put(123,40){${\tiny z}$}
\put(105,35){\vector(0,-1){60}}
\put(90,40){${\tiny z+w}$}
\put(65,35){\vector(0,-1){60}}
\put(20,40){${\tiny z+(n-1)w}$}
\put(138,18){$a$}\put(138,-15){$b$}
\put(80,10){$\cdots$}\put(80,-10){$\cdots$}
\end{picture}
\begin{picture}(50,50)(0,0)
\put(240,5){\vector(1,0){120}}\put(240,5.3){\vector(1,0){120}}
\put(240,5.6){\vector(1,0){120}}\put(240,5.9){\vector(1,0){120}}
\put(325,35){\vector(0,-1){60}}
\put(305,40){${\tiny z+(n-1)w}$}
\put(305,35){\vector(0,-1){60}}
\put(265,35){\vector(0,-1){60}}
\put(257,40){$z$}
\put(338,18){$a$}\put(338,-15){$b$}
\put(280,10){$\cdots$}\put(280,-10){$\cdots$}
\end{picture}
\\ 
\bigskip
\\
$$\begin{array}{cc}
Figure\; 4: The\;definition \;of\quad\quad\quad\quad &\quad\quad\quad\quad
Figure\; 5: The\;definition \;of\\
{\sl\bf L}-matrix \quad\quad\quad & \quad\quad\quad\quad\quad
{\sl\bf L'}-matrix
\end{array}
$$
\\

Under the frame of the theory of interact-round-a-face (IRF), as usual, one can define the  tri-spin operator $\varphi^{(j)}_{a,a+\hat \mu}(z)$ as
$$\varphi^{(j)}_{a,+\hat\mu}(z)=
  \theta\left[\begin{array}{c}
\frac{1}{2}-\frac{j}{n} \\ \frac{1}{2} \end{array}\right](z-nw\bar a_\mu, n\tau),
$$
where $\bar a_\mu=m_\mu-\frac{1}{n}\sum_km_k+\delta_\mu.$ It satisfies the face-vertex corresponding relation (figure 8)
\begin{eqnarray*}
& &\sum_{i_1j_1}r(z_1-z_2)^{i_1 j_1}_{i\ j}
       \varphi^{(i_1)}_{a,a+\hat \mu}(z_1)
       \varphi^{(j_1)}_{a+\hat \mu,a+\hat \mu+\hat \nu}(z_2)\nonumber\\
&=&\sum_kW\left.\left(\begin{array}{cc}
  a & a+\hat k \\ a+\hat\mu & a+\hat\mu+\hat\nu \end{array}\right|  
       z_1-z_2\right)
      \varphi^{(i)}_{a+\hat k,a+\hat \mu+\hat \nu}(z_1)
      \varphi^{(j)}_{a,a+\hat k}(z_2)
\end{eqnarray*}
where  $r$ is the Belavin's $Z_n$ symmetric $R$-matrix and $W$ is the Boltzmann weight of the IRF model $A_{n-1}^{(1)}$. Here we note that in this paper, we express the Boltzmann weight as $R$ in stead of $W$ for convenience.
\pagebreak

\begin{picture}(50,50)(0,0)
\put(40,5){\vector(1,0){120}}
\put(125,35){\vector(0,-1){60}}
\put(105,35){\vector(0,-1){60}}
\put(70,35){\vector(0,-1){60}}
\put(28,8.5){$j'$}\put(152,8.5){$j$}
\put(45,30){$i'_{n-1}$}\put(45,-27){$i_{n-1}$}
\put(80,10){$\cdots$}\put(80,-10){$\cdots$}
\put(95,30){$i'_1$}\put(95,-27){$i_1$}
\put(115,30){$i'_0$}\put(115,-27){$i_0$}
\put(138,18){$a$}\put(138,-15){$b$}
\put(260,5){\vector(1,0){80}}
\put(290,-35){\vector(0,1){80}}
\multiput(290,-32)(0,5){13}{\line(1,1){25}}
\put(300,-10){$a$}\put(300,15){$b$}
\put(270,-5){$j$}\put(250,5){$z$}
\end{picture}
\\[10mm] \mbox{}
$$\begin{array}{cc}
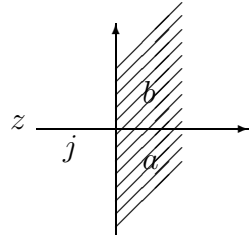
Figure\ 6: The\ definition\ of\quad\quad\quad &\quad\quad\quad
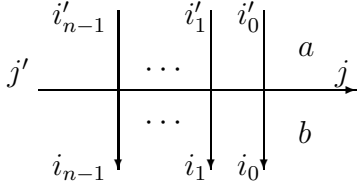
 Figure\ 7: The\ definition\ of \\
X(a|z)^{i_0'i_1'\cdots i_{n-1}'j'}_{i_0i_1\cdots i_{n-1}j} &\quad\quad
\varphi_{a,b}^{(j)}(z)
\end{array}
$$
\\[10mm]
$$
\begin{picture}(50,50)(0,0)
\put(-100,-15.5){\vector(2,1){100}}
\put(-100,15.5){\vector(2,-1){100}}
\put(-30,-50){\vector(0,1){100}}
\multiput(-30,-45)(0,5){17}{\line(1,1){25}}
\put(-27,-35){$a$}
\put(-120,19){$z_1$}\put(-120,-19){$z_2$}
\put(-25,40){$a+\hat\mu+\hat\nu$}
\put(-23,0){$a+{\hat \mu}$}
\end{picture}
=
\begin{picture}(50,50)(0,0)
\put(50,-30.5){\vector(2,1){100}}   
\put(50,30.5){\vector(2,-1){100}}
\put(70,-50){\vector(0,1){100}}
\multiput(70,-45)(0,5){17}{\line(1,1){60}}
\put(110,-25){$a$}
\put(73,-5){$a+\hat k$}
\put(125,0){$a+\hat \mu$}
\put(90,30){$a+\hat\mu+\hat\nu$}
\put(38,33){$z_1$}\put(38,-37){$z_2$}
\end{picture}
$$
\\
\\
$$ Figure\; 8: The\;\; \;face-veter \;\;corresponding \;\;relation.$$

Since the $R$-matrix can be regarded as a special case of ${\bf L}$-matrix, the Eq.(10) is also held for 
\begin{eqnarray}
& &X(a|z)_{i_0\ i_1\ \cdots\ i_{n-1}\ j}^{i'_0\ i'_1\ \cdots\ i'_{n-1}\ j'}
  \nonumber\\
&=&
R(a|z)^{i_0'\ j_1}_{i_0\ j}R(a+\hat i_0'|z+w)^{i_1'\ j_2}_{i_1\ j_1}
R(a+\hat i'_0+\hat i_1'|z+2w)^{i_2'\ j_3}_{i_2\ j_2}\nonumber\\
& &\times\cdots R(a+\hat i'_0+\hat i_1'+\cdots +\hat i_{n-2}'|z+(n-1)w)
   ^{i_{n-1}'\ j'}_{i_{n-1}\ j_{n-1}},
\end{eqnarray}
namely, 
\begin{equation}
P_-(12\cdots n)X(a|z)=P_-(12\cdots n)X(a|z)P_-(12\cdots n),
\end{equation}
giving
\begin{eqnarray}
& &P_-(12\cdots n)^{i_0\ i_1\ \cdots\ i_{n-1}}_{0\ \ 1\ \cdots\ (n-1)}
X(a|z)^{i'_0\ i'_1\ \cdots\ i'_{n-1}\ j'}_{i_0\ i_1\ \cdots\ i_{n-1}\ j}
\varphi^{(0)}_{a,i_0'}(z)\varphi^{(1)}_{a+i_0',i_1'}(z+w)\nonumber\\
& &\times\varphi^{(2)}_{a+i_0'+i_1',i_2'}(z+2w)\cdots
\varphi^{(n-1)}_{a+i_0'+i_1'+\cdots +i_{n-2}',i_{n-1}'}(z+(n-1)w)\nonumber\\
& &\times\varphi^{(i)}_{a+i_0'+i_1'+\cdots+i_{n-1}',j'}(0)\nonumber\\
&=&n!(P_-X(a|z)^{j'}_jP_-)^{0\ 1\ 2\ \cdots\ (n-1)}
       _{0\ 1\ 2\ \cdots\ (n-1)}
P_-(12\cdots n)^{i'_0\ i'_1\ \cdots\ i'_{n-1}}_{0\ 1\ \ \cdots\ (n-1)}
 \nonumber\\ 
& &\times\varphi^{(0)}_{a,i_0'}(z)\varphi^{(1)}_{a+i_0',i_1'}(z+w) \varphi^{(2)}_{a+i_0'+i_1',i_2'}(z+2w)\cdots\nonumber\\
& &\times\varphi^{(n-1)}_{a+i_0'+i_1'+\cdots +i_{n-2}',i_{n-1}'}(z+(n-1)w)
\varphi^{(i)}_{a+i_0'+i_1'+\cdots+i_{n-1}',j'}(0).
\end{eqnarray}
  Obviously, only when $(i_0' i_1' i_2'\cdots i_{n-1}')$ is the permutation of $(012\cdots (n-1))$, the right hand side of Eq.(13) is not zero, otherwise, $P_-(012\cdots (n-1))^{i'_0\ i'_1\ \cdots\ i'_{n-1}}_{0\ 1\ \ \cdots\ (n-1)}
=0$. This implies that only when $j=j'$, the right hand side of Eq.(13) is non-zero. Similar to the Eq. (10), applying $P_-$ on $\varphi$, we obtain from the face-vertex correspondence 
\begin{eqnarray}
&&P_-\varphi\varphi \cdots\varphi
=B^{-1}A\varphi\varphi \cdots\varphi
= B^{-1}\varphi\varphi \cdots\varphi A' \nonumber\\
&=&B^{-1}\varphi\varphi \cdots\varphi A' P_-
=B^{-1}A\varphi\varphi \cdots\varphi P_-
=P_-\varphi\varphi \cdots\varphi P_-,
\end{eqnarray}
where $A'$ is the Cherednik operator constructed by Belavin's $R$-matrix $r$ which also satisfies $r(-w)P=-r(-w)$. The left hand side of  Eq.(13) gives from the face-vertex correspondence 
\begin{eqnarray}
LHS&=&P_-(12\cdots n)^{i_0\ i_1\ \cdots\ i_{n-1}}
        _{0\ \ 1\ \ \cdots\  (n-1)}
\varphi^{(k_0)}_{a+j,i_0}(z)\varphi^{(k_1)}_{a+j+i_0,i_1}(z+w)\nonumber\\
& &\times\varphi^{(k_2)}_{a+j+i_0+i_1,i_2}(z+2w)\cdots
\varphi^{(k_{n-1})}_{a+j+i_0+i_1+\cdots+i_{n-2},i_{n-1}}(z+(n-1)w)
\varphi^{(i')}_{a,j}(0)\nonumber\\
& &\times r(z)^{0\ \ i'_1}_{k_0\ i'}r(z+w)^{1\ \ i'_2}_{k_1\ i'_1}
 r(z+2w)^{2\ \ i'_3}_{k_2\ i'_2}
 \cdots 
r(z+(n-1)w)^{n-1\  i}_{k_{n-1}\ i'_{n-1}}\nonumber\\
&\stackrel{Eq.(14)}{=}&P_-(12\cdots n)^{i_0\ i_1\ \cdots\ i_{n-1}}
        _{0\ \ 1\ \ \cdots\ (n-1)}
\varphi^{(l_0)}_{a+j,i_0}(z)\varphi^{(l_1)}_{a+j+i_0,i_1}(z+w)\nonumber\\
& &\times\varphi^{(l_2)}_{a+j+i_0+i_1,i_2}(z+2w)\cdots
\varphi^{(l_n-1)}_{a+j+i_0+i_1+\cdots+i_{n-2},i_{n-1}}(z+(n-1)w)
\nonumber\\
& &\times P_-(12\cdots (n-1))^{k_0\ k_1\ \cdots\ k_{n-1}}
             _{l_0\ \ l_1\ \ \cdots \ l_(n-1)} r(z)^{0\ \ i'_1}_{k_0\ i'}r(z+w)^{1\ \ i'_2}_{k_1\ i'_1}\nonumber\\
& &\times r(z+2w)^{2\ \ i'_3}_{k_2\ i'_2}
 \cdots 
r(z+(n-1)w)^{n-1\ i}_{k_{n-1}\ i'_{n-1}}\varphi^{(i')}_{a,j}(0).
\end{eqnarray}
One can show  that in the above relation, we can write $P_-rr\cdots r=f(z)\delta^{i}_{i'}P_-$. Thus, comparing the Eq.(13) and Eq.(15), we obtain
\begin{eqnarray*}
& &(P_-X(a|z)^{j}_jP_-)^{012\cdots (n-1)}_{012\cdots(n-1)}\delta^{j'}_j
\varphi^{(i)}_{a+i'_0+i'_1+\cdots+i'_{n-1},j'}(0)
P_-(12\cdots n)^{i'_0\ i'_1\ \cdots\ i'_{n-1}}_{0\ 1\ \ \cdots\ (n-1)}
 \nonumber\\ 
& &\times\varphi^{(0)}_{a,i_0'}(z)\varphi^{(1)}_{a+i_0',i_1'}(z+w) \varphi^{(2)}_{a+i_0'+i_1',i_2'}(z+2w)\cdots
\varphi^{(n-1)}_{a+i_0'+i_1'+\cdots +i_{n-2}',i_{n-1}'}(z+(n-1)w)\\
&\stackrel{Eq.(14)}{=}&P_-(123\cdots n)^{i_0\ i_1\ i_2\ \cdots\ i_{n-1}}
     _{0\ 1\ \ 2\ \ \cdots \ (n-1)}
\varphi^{(0)}_{a+j,i_0}(z)\varphi^{(1)}_{a+j+i_0,i_1}(z+w)
\varphi^{(2)}_{a+j+i_0+i_1,i_2}(z+2w)\nonumber\\
&&\times\cdots\varphi^{(n-1)}
     _{a+j+i_0+i_1+\cdots+i_{n-2},i_{n-1}}(z+(n-1)w)
          f(z)\delta^i_{i'}\varphi^{(i')}_{a,j}(0).
\end{eqnarray*}
From the property of the face wight, we know that
$$
\varphi^{(i)}_{a+i'_0+i'_1+\cdots+i'_{n-1},j'}(0)=
\varphi^{(i)}_{a+\hat 0+\hat 1+\cdots+\hat {n-1},j'}(0)
=\varphi^{(i)}_{a,j'}(0).
$$
Therefore, the above equation gives
\begin{eqnarray}
&& n!(P_-X(a|z)^{j'}_jP_-)^{0\ 1\ 2\ \cdots\ (n-1)}_{0\ 1\ 2\ \cdots\ (n-1)}  
P_-(123\cdots n)^{i'_0\ i'_1\ i'_2\ \cdots\ i'_{n-1}}
      _{0\ \ 1\ 2\ \cdots \ 3}
\varphi^{(0)}_{a,i_0'}(z)\varphi^{(1)}_{a+i_0',i_1'}(z+w)\nonumber\\
& &\times \varphi^{(2)}_{a+i_0'+i_1',i_2'}(z+2w)\cdots
\varphi^{(n-1)}_{a+i_0'+i_1'+\cdots+i_{n-2}',i_{n-1}'}(z+(n-1)w)\nonumber\\
&=&P_-(123\cdots n)^{i_0\ i_1\ i_2 \cdots i_{n-1}}
     _{0\ \ 1\ \  2 \cdots (n-1)}
\varphi^{(0)}_{a+j,i_0}(z)\varphi^{(1)}_{a+j+i_0,i_1}(z+w)
\varphi^{(2)}_{a+j+i_0+i_1,i_2}(z+2w)\nonumber\\
&&\times\cdots\varphi^{(n-1)}
     _{a+j+i_0+i_1+\cdots+i_{n-2},i_{n-1}}(z+(n-1)w)
          f(z)
\end{eqnarray}

By using the definition of the tri-spin operator, we can obtain $\varphi^{(i)}_{a,a+\hat\mu}(z)=\theta^{(i)}(z+nw\bar a_\mu+(n-1)w)$, where $\theta^{(i)}$ is the abbreviation of the $\theta$ function in the definition of $\varphi$. With this notation, we now calculate function $<P_-\varphi P_->$ as follows,
\begin{eqnarray*}
&&<P_-\varphi(a|z)P_->\nonumber\\
 &\equiv & P_-(123\cdots n)^{\mu\nu\lambda\cdots\tau}_{012\cdots(n-1)}
\varphi^{(0)}_{a,\mu}(z)
\varphi^{(1)}_{a+\mu,\nu}(z+w)
\varphi^{(2)}_{a+\mu+\nu,\lambda}(z+2w)\\
& &\times\cdots
\varphi^{(n-1)}_{a+\mu+\nu+\lambda+\cdots,\tau}(z+(n-1)w)\nonumber\\
&=&P_-(1234)^{\mu\nu\lambda\cdots\tau}_{012\cdots(n-1)}
\theta^{(0)}(z+nw\bar a_\mu+(n-1)w)
\theta^{(1)}(z+nw\bar a'_\nu+(n-1)w)\nonumber\\
&&\times\theta^{(2)}(z+nw\bar a''_\lambda+(n-1)w)\cdots
\theta^{(n-1)}(z+nw\bar a^{(n-1)}_\tau+(n-1)w),
\end{eqnarray*}
where $a'=a+\hat\mu,\ a''=a+\hat\mu+\hat\nu,\cdots, \ a^{(n-1)}=a+\hat\mu+\hat\nu+\hat\lambda+\cdots $, and $\hat\mu\ne\hat\nu\ne\hat\lambda\ne\cdots\ne\hat\tau$. Then, we have 
\begin{eqnarray*}
\bar a'_{\nu}&=&m'_{\nu}-\frac{1}{n}\sum_l m'_l+w_{\nu}
               =m'_{\nu}-\frac{1}{n}(\sum_{l\ne\mu} m_l+m_\mu+1)+w_{\nu}\\
             &=&m_{\nu}-\frac{1}{n}(\sum_{l\ne\mu} m_l+m_\mu+1)+w_{\nu}\\
             &=&\bar a_{\nu}-\frac{1}{n}.
\end{eqnarray*}
Similarly, we can also have $a_{\lambda}''=\bar a_{\lambda}-\frac{2}{n},  \cdots,
\bar a^{(n-1)}_{\tau}=\bar a_{\tau}-\frac{n-1}{n}.$
 Let $nz_\mu=z+nw\bar a_\mu+(n-1)$, we then can obtain the final result of $<P_-\varphi P_->$ as
\begin{eqnarray}
& &<P_-\varphi(a|z) P_->\nonumber\\
&=&\mbox{Det}\left[\begin{array}{cccccc}
\theta^{(0)}(nz_0) & \theta^{(1)}(nz_0) & \theta^{(2)}(nz_0) &
\cdots & \theta^{(n-1)}(nz_0)\\
\theta^{(0)}(nz_1) & \theta^{(1)}(nz_1) & \theta^{(2)}(nz_1) &
\cdots & \theta^{(n-1)}(nz_1)\\
\theta^{(0)}(nz_2) & \theta^{(1)}(nz_2) & \theta^{(2)}(nz_2) &
\cdots & \theta^{(n-1)}(nz_2)\\
\vdots & \vdots & \vdots & \ddots & \vdots\\
\theta^{(0)}(nz_{n-1}) & \theta^{(1)}(nz_{n-1}) & \theta^{(2)}(nz_{n-1}) &
\cdots & \theta^{(n-1)}(nz_{n-1}) \end{array}
     \right]\nonumber\\
&=&c\sigma(\sum_i z_i-\frac{n-1}{2})\prod_{i<j}\sigma(z_i-z_j)\nonumber\\
&=&c \sigma(z+w\sum_jw_j+(n-1)(w-\frac{1}{2}))
\prod_{i<j}\sigma(w(\bar a_i-\bar a_j))\nonumber\\
&\equiv& c \sigma(z+w\sum_jw_j+(n-1)(w-\frac{1}{2}))\Delta(a).
\end{eqnarray}
So we obtain from the Eq.(16), 
\begin{equation}
<P_-X(a|z)^{j'}_jP_->
\equiv(P_-X(a|z)^{j'}_jP_-)^{012\cdots (n-1)}_{012\cdots (n-1)}
=\delta^{j'}_j\frac{\Delta(a+\hat j)}{n!\Delta(a)}f(z).
\end{equation}
 By using the DYBR Eq.(2) repeatedly, one can obtain ($j$ and $j'$ are not summed)
$$ 
P_-X(b|z-u)^{j'}_j{\bf L}(^a_b|z)L(^a_b|u)^{j''}_{j'}
=L(^a_b|u)^{j'}_jP_-{\bf L}(^{a+\hat j'}_{b+\hat j}|z)X(a|z-u)^{j''}_j.
$$
Considering the properties of $P_-$, the above equation gives
\begin{eqnarray*}
&&P_-X(b|z-u)^{j'}_jP_-P_-{\bf L}(^a_b|z)P_-L(^a_b|u)^{j''}_{j'}\\
&&=\ \ \ L(^a_b|u)^{j'}_jP_-{\bf L}(^{a+\hat j'}_{b+\hat j}|z)P_-P_-
         X(a|z-u)^{j''}_{j'}P_- \\
&\Rightarrow & 
(P_-X(b|z-u)^{j}_jP_-)\delta^j_{j'}(P_-{\bf L}(^a_b|z)P_-)L(^a_b|u)^{j''}_{j'}\\
&&=\ \ \ L(^a_b|u)^{j}_j(P_-{\bf L}(^{a+\hat j'}_{b+\hat j}|z)P_-)(P_-
         X(a|z-u)^{j'}_{j'}P_-)\delta^{j'}_{j''}\\
&\Rightarrow &
(P_-X(b|z-u)^{j}_jP_-)(P_-{\bf L}(^a_b|z)P_-)L(^a_b|u)^{j''}_j\\
&&=\ \ \ L(^a_b|u)^{j''}_j(P_-{\bf L}(^{a+\hat j''}_{b+\hat j}|z)P_-)(P_-
         X(a|z-u)^{j''}_{j''}P_-).
\end{eqnarray*}
Let $I(^a_b|z)\equiv P_-{\bf L}(^a_b|z)=P_-{\bf L}(^a_b|z)P_-$. Substituting the Eq.(17) into the last relation of above relations, we have
$$\frac{\Delta(b+\hat j)}{\Delta(b)}I(^a_b|z) L(^a_b|u)^{j''}_{j}
=L(^a_b|u)^{j''}_{j}I(^{a+\hat j''}_{b+\hat j}|z)\frac{\Delta(a+\hat j'')}{\Delta(a)}.
$$
Then, we have
\begin{equation}
\left[\frac{\Delta(a)}{\Delta(b)}I(^a_b|z)\right] L(^a_b|u)^{j''}_{j}
=L(^a_b|u)^{j''}_{j}\left[\frac{\Delta(a+\hat j'')}{\Delta(b+\hat j)}
                I(^{a+\hat j''}_{b+\hat j}|z)\right]. 
\end{equation}
Therefore  $(\Delta(a)/\Delta(b))I(^a_b|z)$ is the center in the meaning given by Felder. By using the definition of $P_-$, we can obtain the result of $I$ as follows,
\begin{eqnarray*}
&&I(^a_b|z)=\frac{1}{n!}\sum_P
(-1)^{\left[\mbox{Sign}P(^{0\ 1\ \cdots\
n-1}_{\mu_0\mu_1\cdots\mu_{n-1}})\right]}\\
&&\times\ L(^a_b|z)^0_{\mu_0}L(^{a+\hat 0}_{b+\hat\mu_0}|z+w)^1_{\mu_1}
\cdots L(^{a+\hat 0+\hat 1+\cdots+\hat{n-2}}
_{b+\hat\mu_0+\hat\mu_1+\cdots+\hat\mu_{n-2}}|z+(n-1)w)
^{n-1}_{\mu_{n-1}},
\end{eqnarray*}
and  $P$'s are permutations of integers $0,\ 1,\ \cdots,\ n-1$. 
This agrees with that of Ref.\cite{fe2} for $n=2$.

\section{The coefficient algebra of the $L$-matrix and its center}
~~
For elements of the $L$-matrix which satisfy the Eq.(2), we assume that they take the special forms as follows,
\begin{eqnarray}
 L(^a_b|z)^j_i&=&(^a_b)^j_i\sigma(z+\delta_0+b_i-a_j)F(z),\\
 L(^{a+\hat i'}_{b+\hat i}|z)^{j'}_j
&=&(^{a+\hat i'}_{b+\hat i})^{j'}_j\sigma(z+\delta_0+b'_j-a'_{j'})F(z),
\end{eqnarray}
where $(^a_b)^j_i$ and $(^{a+\hat i'}_{b+\hat i})^{j'}_j$ are the coefficient parts of the elements of $L$-matrix, and they are independent of $z$, $F(z)$ is the function of $z$. Substituting the above equations into the DYBR Eq.(2), we can obtain the following relations (the more detailed derivation is given in Ref.\cite{zsy}) 
\begin{eqnarray}
&&Y^{i'j'}_{i\;i}-Y^{j'i'}_{i\;i}=0, \quad i'\ne j'\\
&&Y^{i'i'}_{i\;j}-Y^{i'i'}_{j\;i}=0, \quad  i\ne j\\
&&\sigma(w)\sigma(a_{i'j'}+b_{ij})Y^{i'j'}_{j\;i}
+\sigma(a_{i'j'})\sigma(b_{ij}-w)Y^{i'j'}_{i\;j}\nonumber\\
&&-\ \sigma(a_{i'j'}+w)\sigma(b_{ij})Y^{j'i'}_{j\;i}=0,
\quad\quad\quad\quad\quad\quad\quad  i\ne i',\ j\ne j'
\end{eqnarray}
where we have used the notation $a_{i'j'}=a_{i'}-a_{j'},\;b_{ij}=b_i-b_j,$ and the definition
\begin{equation}
[^a_b]^{i'}_i[^{a+\hat i'}_{b+\hat i}]^{j'}_j\equiv Y^{i'j'}_{i\;j}
\end{equation}
with 
$$
(^a_b)^{i'}_i\times \prod_{l(\ne i')}\sigma(a_l-a_{i'})\equiv [^a_b]^{i'}_i. 
$$

Eqs.(21)-(23) can
be regarded as the algebraic relations which are satisfied by the operators in the lattice
$(a=\sum_{j=0}^{n-1}m^a_j\hat j, b=\sum_{i=0}^{n-1}m^b_i\hat i)$.
We define a new operator
\begin{equation}
A^{i'}_i\equiv [^a_b]^{i'}_i\Gamma^{i'}_i,
\end{equation}
where
\begin{equation}
 \Gamma^{i'}_if(a,b)=f(a+\hat i',b+\hat i)\Gamma^{i'}_i.
\end{equation}
Namely, we regard the $a,b$ as operators, $\Gamma^{i'}_i$ is not
commutative with the functions of $a,b$. In this way,  we have the
following exchange relations of the operators $\{A^{i'}_i\}$
\begin{eqnarray}
(a)& &A^{i'}_iA^{j'}_i=A^{j'}_iA^{i'}_i,\quad\quad i'\ne j' \nonumber\\
(b)& &\sigma(a_{i'j'}+w)\sigma(b_{ij})A^{j'}_jA^{i'}_i \nonumber\\
   &=&\sigma(a_{i'j'})\sigma(b_{ij}-w)A^{i'}_iA^{j'}_j
     +\sigma(w)\sigma(a_{i'j'}+b_{ij})A^{i'}_jA^{j'}_i,
      \ i\ne i',\  j\ne j'\\
(c)& &A^{i'}_iA^{i'}_j=A^{i'}_jA^{i'}_i. \quad\quad i\ne j\nonumber
\end{eqnarray}
These equations are equivalent relations to the Felder and
Varchenko's elliptic quantum algebra under  special condition. 

For the algebra discussed above, we can prove that it has a set PBW (Poincare-Birkhoff-Witt) base (the precise proof for the existence of PBW base is given in another paper\cite{zsy}). In the following part of this paper, we will discuss the center of this algebra.
 
In the previous section, we have proved that $(\Delta(a)/\Delta(b))I(^a_b|z)$ is the center of the   elliptic quantum group. Here $\Delta(b)\equiv\prod_{i<j}\sigma(w(\bar b_i-\bar b_j))$, the definition of $\Delta(a)$ is similar,  and 
 $I(^a_b|z)$ is given by
\begin{eqnarray*}
&&I(^a_b|z)=\frac{1}{n!}\sum_P
(-1)^{\left[\mbox{Sign}P(^{0\ 1\ \cdots\
n-1}_{\mu_0\mu_1\cdots\mu_{n-1}})\right]}\\
&&\times\ L(^a_b|z)^0_{\mu_0}L(^{a+\hat 0}_{b+\hat\mu_0}|z+w)^1_{\mu_1}
\cdots L(^{a+\hat 0+\hat 1+\cdots+\hat{n-2}}
_{b+\hat\mu_0+\hat\mu_1+\cdots+\hat\mu_{n-2}}|z+(n-1)w)
^{n-1}_{\mu_{n-1}}.
\end{eqnarray*}

In our present case, the Eq.(20) can be rewritten as
$$L(^a_b|z)^{i'}_{i}=\sigma(z+\delta+b_i-a_{i'})A^{i'}_i.$$
So the  quantum determinant of the present algebra can be written as 
\begin{eqnarray}
&&I(^a_b|z)=\frac{1}{n!}\sum_P(-1)^{\left[\mbox{Sign} P(^{0\ 1\ \cdots\
n-1}_{\mu_0\mu_1\cdots\mu_{n-1}})\right]}\nonumber\\ 
&&\times\ \sigma(z+\delta+b_{\mu_0}-a_0)
   \sigma(z+w+\delta+b_{\mu_1}-a_1)\cdots \nonumber\\
&&\times\  \sigma(z+(n-1)w+\delta+b_{\mu_{n-1}}-a_{n-1})
 A^0_{\mu_0}A^1_{\mu_1}\cdots A^{n-1}_{\mu_{n-1}}.
\end{eqnarray}

For each $z$, we can obtain a center element $(\Delta(a)/\Delta(b))I(^a_b|z)$ of algebra (28). We can show that there are only $n$ linearly independent such elements by studying the quasi periodicity of the coefficients\cite{zsy} of Eq.(29)
$$\Phi(z)_{\mu_0\cdots\mu_{n-1}}\equiv \sigma(z+\delta+b_{\mu_0}-a_0)
\cdots \sigma(z+(n-1)w+\delta+b_{\mu_{n-1}}-a_{n-1}).    $$

\section{Conclusion}
~~
The elliptic quantum group is proposed based on the dynamic Yang-Baxter Relation. The generators of the group can be given by the elements of $L$-matrix. we know that the center for a group or algebra play the role that it can commute with all elements of the group or algebra. In the first part of this paper, by using the fusion method, we gave the center for the elliptic quantum group $E_{\tau,\eta}(sl_n)$.  After we  assumed that the elements of $L$-matrix take special forms Eq.(19) and Eq.(20), we obtained a algebra whose algebraic relation  can be regarded as a Yang-Baxter relation. Then, Taking use of the results of the first part, we also obtained the center for the algebra. For this new type of algebra, we can also discuss other properties, such as its base. Actually, we have found that it have a set of PBW base. We will discuss it in our future paper\cite{zsy}.

\end{document}